\documentclass[12pt]{article}
\usepackage{amsfonts}

\usepackage{latexsym}
\usepackage{amssymb}
\usepackage{epsfig}
\usepackage{amsmath}
\usepackage{amsthm}
\usepackage[mathscr]{eucal}
\usepackage{subfigure}
\usepackage{graphicx}

\newtheorem{Lemma}{Lemma}

\newtheorem{Theorem}[Lemma]{Theorem}

\newtheorem{Definition}{Definition}

\renewcommand{\qed}{\hfill{\ \ \rule{2mm}{2mm}} \vspace{0.2in}}

\newcommand{\ind}{1\hspace{-2.3mm}{1}}

\begin{document}

\title{Randomized detection and detection capacity of multidetector networks}
\author{ \textbf{Ghurumuruhan Ganesan}
\thanks{E-Mail: \texttt{gganesan82@gmail.com} } \\
\ \\
New York University, Abu Dhabi}
\date{}

\maketitle

\begin{abstract}



In this paper, we study the following detection problem. There are~\(n\) detectors randomly
placed in the unit square~\(S = \left[-\frac{1}{2},\frac{1}{2}\right]^2\) assigned to detect
the presence of a source located at the origin. Time is divided into slots of unit length
and~\(D_i(t) \in \{0,1\}\) represents the (random)
decision of the~\(i^{th}\) detector in time slot~\(t.\) The location of the source is unknown to the detectors
and the goal is to design
schemes that use the decisions~\(\{D_i(t)\}_{i,t}\) and detect the presence of the source
in as short time as possible. 

We first determine the minimum achievable detection time~\(T_{cap}\) and show the existence of \emph{randomized} detection schemes
that have detection times arbitrarily close to~\(T_{cap}\) for almost all configuration of detectors, provided the
number of detectors~\(n\) is sufficiently large.
We call such schemes as \emph{capacity achieving} and completely characterize all
capacity achieving detection schemes. 

\vspace{0.1in} \noindent \textbf{Key words:} detection capacity, multidetector network.

\end{abstract}

\bigskip



\renewcommand{\theequation}{\thesection.\arabic{equation}}
\setcounter{equation}{0}
\section{Introduction}
\subsection*{Model Description}
Consider \(n\) detectors labelled \(\{1,2,\ldots,n\}\) located in the unit square~\(S = \left[-\frac{1}{2},\frac{1}{2}\right]^2\) and let \(\omega_i, 1 \leq i\ \leq n\) denote the location of the \(i^{th}\) detector. There is also a \emph{source} present at the origin and the location of the source is unknown to the~\(n\) detectors.


The source continuously emits signals and the detectors can therefore sense the presence of the source by receiving and analyzing these signals. We divide time into disjoint slots of unit length and in time slot \(t \geq 1,\) we let~\(D_i(t) \in \{0,1\}\) be the decision of detector~\(i \in \{1,\ldots,n\}\)
regarding the source. Thus \(D_i(t)  = 0\) implies that detector~\(i\) has not detected the source and \(D_i(t) = 1\) implies that user \(i\) has detected the source at time \(t.\)


Let \(r \geq 1\) be a fixed integer. We define~\(\underline{D}_i = (D_i(1),\ldots,D_i(r))\) to be the decision vector of detector~\(i\) in a round of duration~\(r.\) We assume that~\(\underline{D}_i\) consists of independent and identically distributed random variables where
\begin{equation}\label{d_i_t}
\mathbb{P}_{i}(D_i(t) = 1) = p_i = 1 - \mathbb{P}_i(D_i(t) = 0)
\end{equation}
for every \(1 \leq t \leq r,\) where~\(p_i\) denotes the detection probability of the detector~\(i\) and does not depend on the time~\(t.\) We further assume~\(\underline{D}_i\) is independent of~\(\underline{D}_j\) for~\(i \neq j.\) We define the decision vectors~\((\underline{D}_1,\ldots,\underline{D}_n)\) on the probability space~\((\Omega_{dec},{\cal F}_{dec}, \mathbb{P}_{dec})\) where~\(\Omega_{dec} = \{0,1\}^{nr},\) \({\cal F}_{dec}\) is the sigma algegra formed by all subsets of~\(\Omega_{dec}\) and~\(\mathbb{P}_{dec} = \otimes_{i=1}^{n} \mathbb{P}_i.\)


A \((n,r)-\)\emph{detection scheme} is a (deterministic) map \(\pi: \{1,2,3,\ldots,r\} \rightarrow \{1,2,\ldots,n\}.\) In other words, the map~\(\pi\) assigns user~\(\pi(t)\) to detect the channel at time slot \(t.\)

Let \(\underline{X} = (D_{\pi(1)}(1),D_{\pi(2)}(2),\ldots,D_{\pi(r)}(r)) \in \{0,1\}^r\) denote the vector containing the corresponding decisions of the users. Throughout the paper we work only with the vector \(\underline{X}.\) If \(D_{\pi(i)}(i) = 1\) for some \(1 \leq i \leq r,\) we say that the source has been detected and define the corresponding event as~\(A_{det}.\) We also define the detection time random variable~\(T_{det}\) as
\begin{equation}\label{t_det}
T_{det} = \min\{1 \leq  i \leq r : D_{\pi(i)}(i) = 1\}.
\end{equation}
If the event \(A^{c}_{det}\) occurs i.e., the source has not been detected in a round of duration~\(r,\) then we set \(T_{det} = \infty.\)

For a fixed \((n,r)-\)detection scheme~\(\pi\) and a fixed detection probability vector \(\underline{p} = (p_1,\ldots,p_n),\) let \( 1-q_{\pi(j)} = p_{\pi(j)}\) denote the detection probability of user~\(\pi(j)\) at time slot \(1 \leq j \leq r.\) For any fixed \(1 \leq k \leq r\) and fixed \((\pi,\underline{p}),\) we let~\(\mathbb{P}^{(\pi,\underline{p})}\) be the probability measure associated with the decision vectors~\(\{\underline{D}_{i}\}_i\emph{}.\) We then obtain from the model description above that
\begin{equation}\label{prob_pi_p}
\mathbb{P}^{(\pi,\underline{p})}(T_{det} = k)  = \prod_{i=1}^{k-1} q_{\pi(i)} (1-q_{\pi(k)})
\end{equation}
and so
\begin{equation}\label{prob_pi_p_inf}
\mathbb{P}^{(\pi,\underline{p})}(T_{det} < \infty) = \mathbb{P}^{(\pi,\underline{p})}(T_{det} \leq r) = 1 -  \prod_{i=1}^{r} q_{\pi(i)}
\end{equation}
is the probability of the event that the source is located in one round (consisting of~\(r\) time slots). We also define
\begin{equation}\label{exp_pi_p}
\mathbb{E}^{(\pi,\underline{p})}(T_{det}\ind(T_{det} < \infty))  = \sum_{k=1}^{r} k\mathbb{P}^{(\pi,\underline{p})}(T_{det} = k)
\end{equation}
is the expected detection time for a fixed pair~\((\pi,\underline{p}).\)

\subsubsection*{Randomness in configuration}
Suppose now we allow for randomness in the detection probability vector to reflect the randomness in the configuration of the detectors. More precisely, we associate with each detector~\(i,\) a random detection probability~\(P_i \in (0,1)\) taking values in a finite set~\(\Omega_{conf}.\) The random detection probability vector~\(\underline{P} = (P_1,\ldots,P_n)\) has independent and identically distributed components and is defined on the probability space \((\Omega^n_{conf}, {\cal F}_{conf},\mathbb{P}_{conf})\) where \({\cal F}_{conf}\) denotes the collection of all subsets of \(\Omega^n_{conf}.\) Thus the equation~(\ref{d_i_t}) holds for a particular realization \(\underline{p} = (p_1,\ldots,p_n)\) of the random vector~\(\underline{P}.\)

We assume \(\Omega_{conf}\) is finite to avoid measure theoretic complications. In practice this could happen, for example, if there is a (finite) grid of possible locations for placing the detectors.

\subsubsection*{Randomness in Detection Schemes}
We now introduce randomness in the detection scheme and study the configuration quenched (i.e., not averaged with respect to configuration) detection times.
Let \(\Omega_{sch} = \Omega_{sch}(n,r)\) denote the set of all \((n,r)-\)detection schemes. A random \((n,r)-\)detection scheme \(\Pi\) is a random element of \(\Omega_{sch}\) defined on the probability space \((\Omega_{sch}, {\cal F}_{sch}, \mathbb{P}_{sch}).\) Here~\(\mathbb{P}_{sch}\) is a probability distribution on~\(\Omega_{sch}.\) Since the source location is not known to the detectors, we would like our detection scheme to be independent of the detector locations. We therefore assume that the random detection scheme~\(\Pi\) is independent of the random tuple~\((\underline{P},\underline{D}).\)

We define the overall detection process on the probability space\\\((\Omega_{tot}, {\cal F}_{tot}, \mathbb{P}_{tot}),\) where \[\Omega_{tot} = \Omega^n_{conf} \times \Omega_{sch} \times \Omega_{dec}, {\cal F}_{tot} = {\cal F}_{conf} \times {\cal F}_{sch} \times {\cal F}_{dec}\] and~\[\mathbb{P}_{tot} = \mathbb{P}_{conf} \times \mathbb{P}_{sch} \times \mathbb{P}_{dec}.\]


\subsection*{Detection Capacity}
We recall that the randomness in the detection probability vector \(\underline{p}\) is caused by the randomness in the configuration (i.e. location) of the detectors. Since the source location is unknown, we would like that any random placement of the detectors yields reasonably low detection time on an average long as we have enough number of detectors. We therefore have the following definition.
\begin{Definition} We say that detection time of \(s  > 0\) is achievable if for every \(\epsilon,\delta > 0,\) there is a \(N = N(\epsilon,\delta) \geq 1\) so that the following holds for \(n \geq N.\) There is a \(r = r(n) \longrightarrow \infty\) as~\(n \rightarrow \infty\) and a probability distribution \(\mathbb{P}_{sch} = \mathbb{P}_{sch}(n,r)\) such that
\begin{equation}\label{cap_def}
\mathbb{P}_{conf}\left(B(\underline{p},s,\epsilon,\delta)\right) > 1- \epsilon.
\end{equation}
where
\begin{equation}
B(\underline{p},s,\epsilon,\delta) = \{\underline{p} \in \Omega_{conf} : S(\underline{p}) > 1-\epsilon \text{ and } T(\underline{p}) < s + \delta \}. \label{b_def}
\end{equation}
Here for a fixed configuration~\(\underline{p},\) the term
\begin{equation}\label{s_def}
S(\underline{p}) = \sum_{\pi} \mathbb{P}^{(\pi,\underline{p})}(T_{det} < \infty)\mathbb{P}_{sch}(\pi)
\end{equation}
denotes the probability that detection occurs within one round, averaged over all possible detection schemes. Similarly,
\begin{equation}\label{t_def}
T(\underline{p}) = \sum_{\pi} \mathbb{E}^{(\pi,\underline{p})}(T_{det}\ind(T_{det}  < \infty))\mathbb{P}_{sch}(\pi)
\end{equation}
denotes the corresponding averaged detection time for a fixed configuration~\(\underline{p}.\)
\end{Definition}

If~\(\Pi\) is any random~\((n,r)-\)detection scheme with distribution~\(\mathbb{P}_{sch},\) we then say that a detection time of \(s > 0\) is achievable by~\(\Pi\) if the above conditions are satisfied. Roughly speaking, for any random placement of the detectors, the following two conditions must be satisfied with high probability: \((a)\) detection time is finite (i.e., detection happens within one round) and \((b)\) the expected finite detection time is arbitrarily close to~\(s.\)

Define
\begin{equation}\label{t_cap_def}
T_{cap} = \inf\{s > 0: s \text{ is achievable}\}
\end{equation}
to be the detection capacity.
We have the following result regarding the detection capacity.
\begin{Theorem}\label{main_thm}
We have that
\begin{equation}\label{main_est}
T_{cap} = \frac{1}{p_{av}}
\end{equation}
where~
\begin{equation}\label{ave_p2}
p_{av} = (\Omega_{conf})^{-1}\sum_{p \in \Omega_{conf}} p
\end{equation}
is the configuration averaged detection probability. Moreover, a detection scheme~\(\Pi\) achieves a detection time of~\(T_{cap}\) if and only if the following two conditions hold for any fixed integer~\(k \geq 1.\) \\
\((a1)\) We have
\begin{equation}\label{aj_def}
a_k := \mathbb{P}_{sch}\left(\#\{\Pi(1),\ldots,\Pi(k)\} = k\right) \longrightarrow 1
\end{equation}
as~\(n \rightarrow \infty.\)\\
\((a2)\) If~\(\Pi_1\) and~\(\Pi_2\) are two independent detection schemes having the same distribution as~\(\Pi,\) then
\begin{equation}\label{bj_def}
b_k := \mathbb{P}_{sch}\left(\{\Pi_1(1),\ldots,\Pi_1(k)\} \bigcap \{\Pi_2(1),\ldots,\Pi_2(k)\}  = \emptyset \right) \longrightarrow 1
\end{equation}
as~\(n \rightarrow \infty.\)
\end{Theorem}
The result above essentially provides a limit on the detection capability of multidetector networks. This has applications to spectrum sensing in cognitive radio networks, where vacation time is a critical parameter that affects the performance of the network. For more details, we refer to Haykin~(2005), Tandra and Sahai~(2005) and the survey article by Yucek and Arslan~(2009) and references therein. We also refer to Balister et al (2016) for sensing algorithms in a continuum percolation setting.



The paper is organized as follows. In Section~\ref{pf_min_det}, we prove preliminary estimates needed for the proof of Theorem~\ref{main_thm}. In Section~\ref{pf1}, we prove Theorem~\ref{main_thm}.

\renewcommand{\theequation}{\thesection.\arabic{equation}}
\setcounter{equation}{0}
\section{Preliminary estimates}\label{pf_min_det}
We recall from~(\ref{t_def}) that
\begin{equation}\label{t_p_def2}
T(\underline{p}) = \sum_{\pi} \mathbb{E}^{(\pi,\underline{p})}(T_{det}\ind(T_{det} < \infty)) \mathbb{P}_{sch}(\pi)
\end{equation}
is the detection time for a fixed configuration~\(\underline{p},\) averaged over all possible detection schemes.

\subsection*{Mean of~\(T(\underline{p})\)}
We have the following result.
\begin{Lemma}\label{lem1}
Let~\(p_{av}\) be the configuration averaged detection probability as defined in~(\ref{ave_p2}) and let~\(p_{min} = \min\{p : p \in \Omega_{conf}\} > 0\) be the minimum detection probability. For a fixed~\(\epsilon > 0,\) we have that
\begin{equation}\label{e_conf}
\frac{1}{p_{av}} -\epsilon \leq \mathbb{E}_{conf}T(\underline{p}) \leq \frac{1}{p_{av}} + \frac{1}{p_{min}}
\end{equation}
for all~\(n\) large. Also
\begin{equation}\label{e_conf_conv2}
\mathbb{E}_{conf}T(\underline{p}) \longrightarrow \frac{1}{p_{av}}
\end{equation}
as~\(n \rightarrow \infty\) if and only if for each integer~\(k \geq 1,\) the following condition holds:
\begin{equation}\label{cond1}
a_k \longrightarrow 1
\end{equation}
as~\(n \rightarrow \infty.\) Here~\[a_k= a_k(n) := \mathbb{P}_{sch}\left(\#\{\Pi(1),\ldots,\Pi(k)\} = k \right)\] is as defined in~(\ref{aj_def}).
\end{Lemma}
\emph{Proof of Lemma~\ref{lem1}}: From~(\ref{exp_pi_p}), we have that
\begin{eqnarray}
\mathbb{E}^{(\pi,\underline{p})}(T_{det}\ind(T_{det} < \infty))  &=& \sum_{k=1}^{r} k\mathbb{P}^{(\pi,\underline{p})}(T_{det} = k) \nonumber\\
&=&\sum_{k=1}^{r} k\prod_{i=1}^{k-1} q_{\pi(i)} (1-q_{\pi(k)}) \nonumber\\
&=& \sum_{k=1}^{r} k(\alpha_{k-1} - \alpha_{k}) \nonumber\\
&=& \sum_{j=0}^{r-1} \alpha_j - r\alpha_{r}. \label{e_pi_p_est}
\end{eqnarray}
where~\(\alpha_0 = 1\) and for~\(k \geq 1,\) we have
\begin{equation}\label{alpha_k_def}
\alpha_k  = \alpha_k(\pi,\underline{p}) := \prod_{i=1}^{k} q_{\pi(i)}.
\end{equation}

We have from the definition that~\(\alpha_j\) depends on both the configuration~\(\underline{p}\) and the detection scheme~\(\pi.\) From~(\ref{e_pi_p_est}), we have that the configuration averaged detection time is
\begin{equation}\label{conf_ave_det}
T(\underline{p}) = \sum_{j=0}^{r-1} \mathbb{E}_{sch}(\alpha_j) - r \mathbb{E}_{sch}(\alpha_r)
\end{equation}
and so
\begin{eqnarray}
\mathbb{E}_{conf}T(\underline{p}) &=& \sum_{j=0}^{r-1} \mathbb{E}_{conf}\mathbb{E}_{sch}(\alpha_j) - r \mathbb{E}_{conf}\mathbb{E}_{sch}(\alpha_r) \nonumber\\
&=& \sum_{j=0}^{r-1} \mathbb{E}_{sch}\mathbb{E}_{conf}(\alpha_j) - r \mathbb{E}_{sch}\mathbb{E}_{conf}(\alpha_r). \label{ov_ave}
\end{eqnarray}

For a fixed detection scheme~\(\pi,\) we first estimate~\(\mathbb{E}_{conf}(\alpha_j).\) For~\(j =1,\) we have that~\[\mathbb{E}_{conf}(\alpha_1) = 1-\mathbb{E}_{conf}(p_{\pi(1)}) =  1-p_{av}\] for any detection scheme~\(\pi,\) where~\(p_{av}\) is the detection probability averaged over all possible configurations as defined in~(\ref{ave_p2}). For a fixed integer~\(j \geq 2\) and a fixed detection scheme~\(\pi,\) we have the following estimates for the configuration averaged value of~\(\alpha_j.\) We have
\begin{equation}\label{conf_ave_alphaj2}
(1-p_{av})^{j} \leq \mathbb{E}_{conf}(\alpha_j) \leq (1-p_{min})^{j}
\end{equation}
and
\begin{equation}\label{conf_ave_alphaj}
c_j\ind(W_j^c) \leq \mathbb{E}_{conf}(\alpha_j) - (1-p_{av})^{j} \leq d_j\ind(W_j^c)
\end{equation}
for all~\(n \geq 1,\) where
\begin{equation}\label{uj_def}
0 < d_j = (1-p_{min})^{j} - (1-p_{av})^{j}  \leq (1-p_{min})^{j}
\end{equation}
and
\begin{equation}\label{cj_def}
c_j = \min_{2 \leq i \leq j} \frac{\mathbb{E}_{conf}(1-p)^{i} - (1-p_{av})^{i}}{(1-p_{av})^{i}} > 0.
\end{equation}
Here~\(W_j = \{\#\{\pi(1),\ldots,\pi(j)\} = j\}\) is the event that first~\(j\) values of~\(\pi\) are all distinct. The term~\(p_{av}\) is the configuration averaged detection probability as defined in~(\ref{ave_p2}) and~\(p_{min} = \min\{p: p \in \Omega_{conf}\} > 0\) is the minimum detection probability. The estimate~(\ref{conf_ave_alphaj}) is slightly more stronger than~(\ref{conf_ave_alphaj2}) and from~(\ref{conf_ave_alphaj}), we obtain that the term~\(\mathbb{E}_{conf}(\alpha_j) = (1-p_{av})^{j}\) if and only if the event~\(W_j\) occurs.  \\\\
\emph{Proof of~(\ref{conf_ave_alphaj2}) and~(\ref{conf_ave_alphaj})}:
Suppose that~\[\{\pi(1),\ldots,\pi(j)\} = \{i_1.x_1,\ldots,i_w.x_w\},\] where~\(\{x_1,\ldots,x_w\}\) are the distinct elements in~\(\{\pi(1),\ldots,\pi(j)\}\) and~\(i_k\) denotes the multiplicity of~\(x_k\) for~\(1 \leq k \leq w,\) satisfying
\begin{equation}\label{mult_ak}
\sum_{k=1}^{w} i_k = j.
\end{equation}
We recall that~\(q_{\pi(i)} = 1-p_{\pi(i)}\) and~\(p_{\pi(i)}\) is the detection probability for detector~\(\pi(i).\) Using the~\(\mathbb{P}_{conf}-\)independence of the detection probabilities\\\(p_{x_1},\ldots,p_{x_w},\) we then have
\begin{eqnarray}
\mathbb{E}_{conf}(\alpha_j) &=& \prod_{k=1}^{w}\mathbb{E}_{conf}(1-p_{x_k})^{i_k} \label{eq_later_nd}\\
&\geq& \prod_{k=1}^{w}\left(\mathbb{E}_{conf}(1-p_{x_k})\right)^{i_k} \label{mid_exp}\\
&=& \prod_{k=1}^{w}(1-p_{av})^{i_k} \nonumber\\
&=& (1-p_{av})^{j} \label{fin_ealp}
\end{eqnarray}
where the final estimate follows from~(\ref{mult_ak}). In the middle step~(\ref{mid_exp}), we use the estimate~\(\mathbb{E}X^{q} \geq (\mathbb{E}X)^{q}\) for any positive random variable~\(X\) and integer~\(q \geq 1.\) Moreover, equality occurs in~(\ref{mid_exp}) if and only if~\(i_k = 1\) for each~\(1 \leq k \leq w.\) This proves the lower bound in~(\ref{conf_ave_alphaj2}) and the equality in~(\ref{conf_ave_alphaj}) if~\(W_j\) occurs.

Suppose now that~\(W_j^c\) occurs. This means that~\(i_k \geq 2\) for some~\(1 \leq k \leq w.\) Suppose~\(i_1 \geq 2.\) Arguing as in~(\ref{fin_ealp}), we get from~(\ref{eq_later_nd}) that
\begin{eqnarray}
\mathbb{E}_{conf}(\alpha_j) &=& \mathbb{E}_{conf}(1-p_{x_1})^{i_1} \prod_{k=2}^{w}\mathbb{E}_{conf}(1-p_{x_k})^{i_k} \nonumber\\
&\geq& \mathbb{E}_{conf}(1-p_{x_1})^{i_1} \prod_{k=2}^{w} (1-\mathbb{E}_{conf}(p_{x_k}))^{i_k} \label{mid_exp2}\\
&=& \mathbb{E}_{conf}(1-p_{x_1})^{i_1} \prod_{k=2}^{w} (1-p_{av})^{i_k} \nonumber\\
&=& \Delta(i_1) \mathbb{E}_{conf}(1-p_{av})^{j}\nonumber
\end{eqnarray}
where~\[\Delta(i_1) = \frac{\mathbb{E}_{conf}(1-p_{x_1})^{i_1}}{(\mathbb{E}_{conf}(1-p_{x_1}))^{i_1}} \geq 1+c_j\] and~\(c_j > 0\) is as defined in~(\ref{cj_def}). This proves the lower bound in~(\ref{conf_ave_alphaj}).

The upper bound in~(\ref{conf_ave_alphaj2}) and~(\ref{conf_ave_alphaj}) follows from~(\ref{eq_later_nd}) and~(\ref{mult_ak}) along with the fact that~\(p_{x_k} \geq p_{min}\) for all~\(1 \leq k \leq w.\)~\(\qed\)

Substituting the bounds for~(\ref{conf_ave_alphaj2}) into~(\ref{ov_ave}), we get
\begin{eqnarray}
\mathbb{E}_{conf}T(\underline{p}) &\geq& \sum_{j=0}^{r-1} (1-p_{av})^{j} - r (1-p_{min})^{r}\\
&=& \frac{1-(1-p_{av})^{r}}{p_{av}} - r(1-p_{min})^{r} \nonumber\\
&\geq& \frac{1}{p_{av}} - \epsilon \label{low_bd_lem1}
\end{eqnarray}
for all~\(n\) large, provided~\(r = r(\epsilon) \geq 1\) is large and fixed. This proves the lower bound in~(\ref{e_conf}).

For the rest, we argue as follows. Using the upper bound bound in~(\ref{conf_ave_alphaj}) in~(\ref{conf_ave_det}) we have
\begin{equation}\label{ov_ave2}
\mathbb{E}_{conf}T(\underline{p}) \leq \sum_{j=0}^{r-1} (1-p_{av})^{j}  + \sum_{ j=0}^{r-1}(1-p_{min})^{j} (1-a_j)
\end{equation}
and using the lower bound in~(\ref{conf_ave_alphaj}) and upper bound in~(\ref{conf_ave_alphaj2}) in~(\ref{conf_ave_det}), we have
\begin{equation}\label{ov_ave3}
\mathbb{E}_{conf}T(\underline{p}) \geq \sum_{j=0}^{r-1} (1-p_{av})^{j}  + \sum_{ j=0}^{r-1}c_j (1-a_j) - r(1-p_{min})^{r-1}
\end{equation}
where the sequence~\[a_j = a_j(n) := \mathbb{P}_{sch}\left(\#\{\Pi_1(1),\ldots,\Pi_1(j)\} = j\right) \leq 1\] is as defined in~(\ref{aj_def}).

From~(\ref{ov_ave2}) and the fact that~\(1-a_j \leq 1,\) we obtain the upper bound in~(\ref{e_conf}). We now prove~(\ref{e_conf_conv2}). Suppose now that~(\ref{cond1}) holds so that~\(a_j \longrightarrow 1\) as~\(n \rightarrow \infty\) for any fixed integer~\(j \geq 1.\) Fixing integer~\(r \geq 1\) large to be determined later, we have that~\(a_j \geq 1-\epsilon\) for all~\(1 \leq j \leq r\) and for all~\(n \geq N(r,\epsilon) \geq 1.\) Using this in~(\ref{ov_ave2}), we have
\begin{eqnarray}
\mathbb{E}_{conf}T(\underline{p}) &\leq& \sum_{j=0}^{r-1} (1-p_{av})^{j}  + \epsilon \sum_{j=0}^{r-1}(1-p_{min})^{j} \nonumber\\
&\leq& \frac{1}{p_{av}} + \frac{\epsilon}{p_{min}} \label{up_bd_lem1}
\end{eqnarray}
for all~\(n\) large. Since~\(\epsilon  >0\) is arbitrary, we obtain~(\ref{e_conf_conv2}) from~(\ref{low_bd_lem1}) and~(\ref{up_bd_lem1}).

Suppose now that~(\ref{cond1}) does not hold so that there is an integer~\(r_0 \geq 1\) and a number~\(\epsilon_0 > 0\) and a sequence~\(\{n_k\}\) such that~\(a_{r_0} = a_{r_0}(n_k) \leq 1-\epsilon_0\) for all~\(k\) large.
Using this in~(\ref{ov_ave3}), we then have
\begin{eqnarray}
\mathbb{E}_{conf}T(\underline{p}) &\geq& \sum_{j=0}^{r-1} (1-p_{av})^{j}  + \sum_{j=0}^{r-1}c_j (1-a_j)  - r(1-p_{min})^{r-1}\label{e_conf_lower}.
\end{eqnarray}

For~\(\epsilon  >0\) small, we have that \[\sum_{j=0}^{r-1} (1-p_{av})^{j} - r(1-p_{min})^{r-1} \geq \frac{1}{p_{av}} - \epsilon\] for all~\(n\) large provided~\(r = r(\epsilon) \geq 1\) is large. Similarly \[\sum_{j=0}^{r-1} c_j(1-a_j) \geq c_{r_0} (1-a_{r_0}) \geq c_{r_0}\epsilon_0 > 0\] by choice of~\(r_0.\) This implies that
\[\mathbb{E}_{conf}T(\underline{p}) \geq \frac{1}{p_{av}} -\epsilon + c_{r_0}\epsilon_0.\] Since~\(\epsilon  >0\) is arbitrary,~(\ref{e_conf_conv2}) cannot hold.~\(\qed\)

\subsection*{Variance of~\(T(\underline{p})\)}
From~(\ref{conf_ave_det}), we have for a fixed integer~\(k \geq 1\) that
\begin{equation}\label{t_split}
T(\underline{p}) = 1 + \sum_{j=1}^{k} T_j + R_k
\end{equation}
where
\begin{equation}\label{t_j_def}
T_j =  T_j(\underline{p}) = \sum_{\underline{i}} Q(\underline{i}) \beta(\underline{i})
\end{equation}
and~\( \underline{i}  = (i_1,\ldots,i_j) \in \{1,\ldots,n\}^{j}\) is~\(j-\)tuple. 
For a fixed~\(\underline{i}  = (i_1,\ldots,i_j),\) the term
\begin{equation}\label{q_def}
Q(\underline{i}) = \prod_{1 \leq l \leq j} q_{i_l}
\end{equation}
and
\begin{equation}\label{beta_def}
\beta(\underline{i}) = \mathbb{P}_{sch}\left(\Pi(1) = i_1,\ldots,\Pi(j) = i_j\right)
\end{equation}
Similarly the term
\begin{equation}\label{r_k_def}
R_k = \sum_{j=k+1}^{r-1}T_j - rT_r.
\end{equation}

We have the following estimates regarding the~\(Q-\)terms.
\begin{Lemma}\label{var_q}
Fix~\(j_1,j_2 \geq 1\) and~\(\underline{i}_1 \in \{1,\ldots,n\}^{j_1}\) and~\(\underline{i}_2 \in \{1,\ldots,n\}^{j_2}.\)
Let
\begin{equation}\label{q_corr}
\delta(\underline{i}_1,\underline{i}_2) := \mathbb{E}_{conf}Q(\underline{i}_1)Q(\underline{i}_2) - \mathbb{E}_{conf}Q(\underline{i}_1)\mathbb{E}_{conf}Q(\underline{i}_2).
\end{equation}
We have that~\(\delta(\underline{i}_1,\underline{i}_2) = 0\) if and only if~\(\{\underline{i}_1\} \cap \{\underline{i}_2\} = \emptyset;\) i.e., the tuples~\(\underline{i}_1\) and~\(\underline{i}_2\) have no entries in common. Also
\begin{equation}\label{del_est}
\delta(\underline{i}_1,\underline{i}_2) \geq e(j_1,j_2)\ind\left(\{\underline{i}_1\} \cap \{\underline{i}_2\} \neq \emptyset\right)
\end{equation}
for some constant~\(e(j_1,j_2) > 0.\) Moreover
\begin{equation}\label{cross_est1}
\mathbb{E}_{conf}Q(\underline{i}_1) \leq \mathbb{E}_{conf}(1-p)^{j_1}.
\end{equation}
and
\begin{equation}\label{cross_est2}
\mathbb{E}_{conf}Q(\underline{i}_1) Q(\underline{i}_2) \leq \mathbb{E}_{conf}(1-p)^{j_1+j_2}.
\end{equation}
\end{Lemma}
\emph{Proof of Lemma~\ref{var_q}}: Let~\(\{\underline{i}_1\}\) represent the set of indices present in the\\\(j-\)tuple \(\underline{i}_1.\) We have
\[Q(\underline{i}_1) = Q_1Q_{12} \text{ and } Q(\underline{i}_2) = Q_2 Q_{12}\] where~\(Q_1\) represents the product corresponding to indices in~\(\{\underline{i}_1\}\) but not in~\(\{\underline{i}_2\}\) and~\(Q_{12}\) represents the product corresponding to indices present in both~\(\{\underline{i}_1\}\) and~\(\{\underline{i}_2\}.\) Using the~\(\mathbb{P}_{conf}-\)independence of the terms~\(Q_1,Q_{12}\) and~\(Q_2,\) we have
\begin{equation}\label{q_12}
\mathbb{E}_{conf}Q(\underline{i}_1) Q(\underline{i}_2) = \mathbb{E}_{conf}Q_1\mathbb{E}_{conf}Q^2_{12} \mathbb{E}_{conf}Q_2.
\end{equation}
Similarly we have
\begin{equation}\label{q_1_q_2}
\mathbb{E}_{conf}Q(\underline{i}_1) \mathbb{E}_{conf}Q(\underline{i}_2) = \mathbb{E}_{conf}Q_1\left(\mathbb{E}_{conf}Q_{12} \right)^{2} \mathbb{E}_{conf}Q_2.
\end{equation}
If~\(Q_{12} = 1;\) i.e., the tuples~\(\underline{i}_1\) and~\(\underline{i}_2\) have no entries in common, then the term~\(\delta(\underline{i}_1,\underline{i}_2) = 0.\)

If~\(\{\underline{i}_1\} \cap \{\underline{i}_2\} \neq \emptyset,\) then~\(Q_{12} < 1\) strictly and we have from~(\ref{q_12}) and~(\ref{q_1_q_2}) that
\[\delta(\underline{i}_1,\underline{i}_2)  = \mathbb{E}_{conf}(Q_1)var_{conf}(Q_{12})\mathbb{E}_{conf}(Q_2).\] Taking minimum over all possible choices of~\(Q_1,Q_{12}\) and~\(Q_3\) we obtain the lower bound in~(\ref{del_est}). 

We prove~(\ref{cross_est1}) and the proof for~(\ref{cross_est2}) is analogous. Suppose~\(\{\underline{i}_1\} = \{w_1.c_1,\ldots,w_q.c_q\}\) where~\(\{c_i\}\) are the indices present in the tuple~\(\underline{i}_1\) with~\(\{w_i\}\) representing the corresponding multiplicities so that~
\begin{equation}\label{sum_mult}
\sum_{l=1}^{q} w_l = j.
\end{equation}
We then have
\begin{equation}\label{inter1}
\mathbb{E}_{conf}Q(\underline{i}_1) = \prod_{l=1}^{q} \mathbb{E}_{conf}(1-p)^{w_l}.
\end{equation}
Using~\(\left(\mathbb{E}X^{s}\right)^{\frac{1}{s}} \leq \left(\mathbb{E}X^{t}\right)^{\frac{1}{t}}\) for~\(s \leq t\) and a positive random variable~\(X,\) we obtain that
\begin{equation}\label{inter2}
\mathbb{E}_{conf}(1-p)^{w_l} \leq \left(\mathbb{E}_{conf}(1-p)^{j}\right)^{\frac{w_l}{j}}
\end{equation}
and so the final term in~(\ref{inter1}) is at most
\[\prod_{l=1}^{q} \left(\mathbb{E}_{conf}(1-p)^{j}\right)^{\frac{w_l}{j}} = \mathbb{E}_{conf}(1-p)^{j}.\] The final estimate follows from~(\ref{sum_mult}).~\(\qed\)

\begin{Lemma}\label{var_tj}
Fix~\(j_1,j_2 \geq 1.\) We have
\begin{equation}\label{pos_corr}
cov_{conf}\left(T_{j_1},T_{j_2}\right) := \mathbb{E}_{conf}(T_{j_1}T_{j_2}) - \mathbb{E}_{conf}(T_{j_1})\mathbb{E}_{conf}(T_{j_2}) \geq 0.
\end{equation}
Also
\begin{equation}\label{cross_t2}
\mathbb{E}_{conf}(T_{j_1}T_{j_2}) \leq \mathbb{E}_{conf}(1-p)^{j_1+j_2}
\end{equation}
and
\begin{equation}\label{var_est_tj}
e(j_1,j_1) (1-b_{j_1}) \leq var_{conf}(T_{j_1}) \leq 1-b_{j_1}
\end{equation}
where~\(e(j_1,j_1) > 0\) is the constant defined in~(\ref{del_est}) and~\(b_{j_1}\) is the constant as defined in~(\ref{bj_def}).

Fix~\(\epsilon > 0\) and~\(1 \leq k \leq n.\) We have
\begin{equation}\label{var_tp_est_low}
var_{conf}(T(\underline{p})) \geq \sum_{j=1}^{k} e(j,j) (1-b_j) - \epsilon
\end{equation}
for all~\(n \geq N(k,\epsilon) \geq 1\) large.
If~\(k  = k(\epsilon) \geq 1\) is large, we also have
\begin{equation}\label{var_tp_est_up}
var_{conf}(T(\underline{p})) \leq (k+1)\sum_{j=1}^{k} (1-b_j) + \epsilon
\end{equation}
for all~\(n \geq N(k,\epsilon) \geq 1\) large.
\end{Lemma}

\emph{Proof of Lemma~\ref{var_tj}}: The estimate~(\ref{pos_corr}) follows from~(\ref{q_corr}) in Lemma~\ref{var_q} since we have from~(\ref{t_j_def}) that
\begin{equation}
\mathbb{E}_{conf}(T_{j_1}T_{j_2}) - \mathbb{E}_{conf}(T_{j_1})\mathbb{E}_{conf}(T_{j_2})
= \sum_{\underline{i}_1}\sum_{\underline{i}_2} \delta(\underline{i}_1,\underline{i}_2)\beta_n(\underline{i}_1)\beta_n(\underline{i}_2)
\end{equation}
where~\(\delta(.,.) \geq 0\) is as defined in~(\ref{q_corr}) and~\(\underline{i}_1 \in \{1,\ldots,n\}^{j_1},\underline{i}_2 \in \{1,\ldots,n\}^{j_2}.\) Similarly~(\ref{cross_t2}) follows from~(\ref{cross_est2}).

To estimate the upper bound for the variance of~\(T_j\) we proceed as follows.
We have from~(\ref{t_j_def}) that
\[var_{conf}(T_j) = \sum_{\underline{i}_1} \sum_{\underline{i}_2} \delta(\underline{i}_1,\underline{i}_2) \beta(\underline{i}_1)\beta(\underline{i}_2)\]
where~\(\delta(.,.)\) is as defined in~(\ref{q_corr}). From~(\ref{q_corr}) we also have that~\(\delta(\underline{i}_1,\underline{i}_2)  = 0\) if and only if~\(\underline{i}_1\) and~\(\underline{i}_2\) do not have any entries in common; i.e., the sets~\(\{\underline{i}_1\} \cap \{\underline{i}_2\} = \emptyset.\) So
\begin{eqnarray}
var_{conf}(T_j) &=& \sum_{\underline{i}_1,\underline{i}_2 : \{\underline{i}_1\} \cap \{\underline{i}_2\} \neq \emptyset} \delta(\underline{i}_1,\underline{i}_2) \beta(\underline{i}_1)\beta(\underline{i}_2)\nonumber\\
&\leq& \sum_{\underline{i}_1,\underline{i}_2 : \{\underline{i}_1\} \cap \{\underline{i}_2\} \neq \emptyset} \beta(\underline{i}_1)\beta(\underline{i}_2)\nonumber\\
&=& 1-b_j \nonumber
\end{eqnarray}
where the middle estimate follows since~\(\delta(.,.) \leq 1.\) The lower bound similarly follows from the lower bound for~\(\delta(.,.)\) in~(\ref{del_est}).

To prove the lower bound in~(\ref{var_tp_est_low}), we argue as follows. From~(\ref{t_split}) we have
\[T(\underline{p}) = 1 + \sum_{j=1}^{r-1} T_j - rT_r\] and so if~\(r \geq k,\) we have
\begin{eqnarray}
var_{conf}(T(\underline{p})) &=& \sum_{j=1}^{r-1} var_{conf}(T_j) +r^2 var_{conf}(T_r)  - 2r\sum_{j=1}^{r} cov_{conf}(T_j,T_r) \nonumber\\
&&\;\;\;\;\;+ \sum_{1 \leq j_1 \neq j_2 \leq r-1} cov(T_{j_1},T_{j_2}) \nonumber\\
&\geq& \sum_{j=1}^{k} var_{conf}(T_j)  - 2r\sum_{j=1}^{r} cov_{conf}(T_j,T_r) \label{cov_1}\\
&\geq& \sum_{j=1}^{k} var_{conf}(T_j)  - 2r\sum_{j=1}^{r} \mathbb{E}_{conf}T_jT_r. \label{cov_2}
\end{eqnarray}
The inequality~(\ref{cov_1}) follows using~(\ref{pos_corr}) and the estimate~(\ref{cov_2}) follows from the definition of covariance in~(\ref{pos_corr}).

Using~(\ref{cross_t2}) we have
\[\mathbb{E}_{conf}T_jT_r \leq \mathbb{E}_{conf}(1-p)^{j+r} \leq \mathbb{E}_{conf}(1-p)^{r}\]
and using the above in~(\ref{cov_2}), we have
\begin{eqnarray}
var_{conf}(T(\underline{p})) &\geq& \sum_{j=1}^{k} var_{conf}(T_j)  - 2r^2\mathbb{E}_{conf}(1-p)^{r} \nonumber\\
&\geq& \sum_{j=1}^{k} var_{conf}(T_j)  - \epsilon  \nonumber\\
&\geq& \sum_{j=1}^{k} e(j,j)(1-b_j) - \epsilon \label{cov_5}
\end{eqnarray}
provided~\(r = r(\epsilon) \geq 1\) is large. The final estimate follows from the lower bound for the variance of~\(T_j\) in~(\ref{var_est_tj}). This proves the lower bound in~(\ref{var_tp_est_low}).

To prove the upper bound in~(\ref{var_tp_est_up}), we argue as follows. Using~(\ref{t_split}) and the identity~\((\sum_{1 \leq i \leq k}v_i)^2 \leq k\sum_{1 \leq i \leq k} v_i^2,\) we have
\begin{equation}\label{var_tp_est_up2}
var_{conf}(T(\underline{p})) \leq (k+1)\sum_{j=1}   ^{k} var_{conf}(T_j) + (k+1)var_{conf}(R_k).
\end{equation}
We have
\begin{eqnarray}
var_{conf}(R_k) &\leq& \mathbb{E}_{conf}R_k^2 \nonumber\\
&=& \mathbb{E}_{conf}\left(\sum_{j=k+1}^{r-1}T_j\right)^{2} \nonumber\\
&=& \sum_{j_1 = k+1}^{r-1}\sum_{j_2 = k+1}^{r-1}\mathbb{E}_{conf}T_{j_1}T_{j_2} \nonumber\\
&\leq& \sum_{j_1 = k+1}^{r-1}\sum_{j_2 = k+1}^{r-1}\mathbb{E}_{conf}(1-p)^{j_1+j_2} \label{r_k_est}
\end{eqnarray}
where the final estimate follows using~(\ref{cross_est2}).

Using the geometric summation formula we have that the final term in~(\ref{r_k_est}) is
\[\mathbb{E}_{conf} \sum_{j_1 = k+1}^{r-1}\sum_{j_2 = k+1}^{r-1}(1-p)^{j_1+j_2}  \leq \mathbb{E}_{conf} \frac{(1-p)^{2k+2}}{p^2} \leq \frac{1}{p_{min}^2}\mathbb{E}(1-p)^{2k+2}.\] Substituting the above into~(\ref{var_tp_est_up2}), we have
\begin{eqnarray}
var_{conf}(T(\underline{p})) &\leq& (k+1)\sum_{j=1}^{k} var_{conf}(T_j) + (k+1)\frac{1}{p_{min}^2}\mathbb{E}(1-p)^{2k+2} \nonumber\\
&\leq& (k+1)\sum_{j=1}^{k} var_{conf}(T_j) + \epsilon \nonumber\\
&\leq& (k+1)\sum_{j=1}^{k} (1-b_j) + \epsilon \label{est_up3}
\end{eqnarray}
for all~\(n\) large, provided~\(k = k(\epsilon) \geq 1\) is large. The final estimate in~(\ref{est_up3}) follows from the upper bound for the variance of~\(T_j\) in~(\ref{var_est_tj}).~\(\qed\)

The following is the main result of this subsection.
\begin{Lemma}\label{lem2}
Let~\(\Pi\) be a random detection scheme with distribution~\(\mathbb{P}_{sch}.\) The following conditions are equivalent.\\
\((i)\) The term
\begin{equation}\label{cond_i}
var_{conf}(T(\underline{p})) \longrightarrow 0
\end{equation}
as~\( n \rightarrow \infty.\)\\
\((ii)\) The term
\begin{equation}\label{cond_ii}
T(\underline{p}) - \mathbb{E}_{conf}T(\underline{p}) \longrightarrow 0
\end{equation}
in probability, as~\(n \rightarrow \infty.\)\\
\((iii)\) For every fixed~\(k \geq 1,\) we have \(b_k  = b_k(n) \longrightarrow 1\)
as~\(n \rightarrow \infty.\) Here~\(b_k\) is as defined in~(\ref{bj_def}).
\end{Lemma}
\emph{Proof of Lemma~\ref{lem2}}: From~(\ref{var_tp_est_up}), we have that~\(T(\underline{p}) = T_n(\underline{p})\) is a sequence of uniformly integrable (u.i.) random variables and so condition~\((i)\) is equivalent to condition~\((ii).\) Again using the upper bound in~(\ref{var_tp_est_up}), we have that if condition~\((iii)\) holds, then condition~\((i)\) holds.

Suppose now that condition~\((iii)\) does not hold. There exists an integer~\(k_0 \geq 1, \epsilon_0 > 0\) and a sequence~\(\{n_j\}\) such that~\(b_{k_0} = b_{k_0}(n_j) \leq 1-\epsilon_0\) for all~\(j\) large. From the lower bound in~(\ref{var_tp_est_low}), we then have that condition~\((i)\) also does not hold.~\(\qed\)

We need the following Lemma for future use.
\begin{Lemma}\label{lem_yn} Let~\(\{Y_n\}\) be a set of random variables with~\(\mu_n = \mathbb{E} Y_n\) and\\\(\sup_n \mathbb{E}Y_n^2 <\infty.\) Suppose for every~\(\epsilon > 0\) we have
\begin{equation}\label{cond_yn}
\mathbb{P}\left(Y_n < \mu_n + \epsilon\right) \longrightarrow 1
\end{equation}
as~\(n \rightarrow \infty.\) We then have for every~\(\epsilon  >0\) that
\begin{equation}\label{res_yn}
\mathbb{P}\left(Y_n > \mu_n - \epsilon\right) \longrightarrow 1
\end{equation}
as~\(n \rightarrow \infty.\)
\end{Lemma}
\emph{Proof of Lemma~\ref{lem_yn}}: Suppose~(\ref{res_yn}) does not hold. There exists~\(\epsilon_0,\delta_0 >0\) and a subsequence~\(\{n_k\}\) such that
\begin{equation}\label{del_0_def}
\mathbb{P}\left(Y_{n_k} > \mu_{n_k} - \epsilon_0\right) \leq 1-\delta_0
\end{equation}
for all~\(k\) large. Letting~\(F_{n_k} = \{Y_{n_k} > \mu_{n_k} - \epsilon_0\},\) we then have
\begin{eqnarray}
\mathbb{E}Y_{n_k} &=& \mathbb{E}Y_{n_k}\ind(F_{n_k}) + \mathbb{E}Y_{n_k}\ind(F^{c}_{n_k}) \nonumber\\
&\leq& \mathbb{E}Y_{n_k}\ind(F_{n_k}) + (\mu_{n_k} - \epsilon_0)\mathbb{P}(F_{n_k}^c)\label{up1}
\end{eqnarray}
We evaluate the first term in~(\ref{up1}) as follows.
Fix~\(\epsilon  >0\) and let~\(G_{n_k} = \{Y_{n_k} < \mu_{n_k} + \epsilon\}.\) We have that
\begin{eqnarray}
\mathbb{E}Y_{n_k}\ind(F_{n_k}) &=& \mathbb{E}Y_{n_k} \ind(F_{n_k} \cap G_{n_k}) + \mathbb{E}Y_{n_k} \ind(F_{n_k} \cap G_{n_k}) \nonumber\\
&\leq& (\mu_{n_k} + \epsilon)\mathbb{P}(F_{n_k} \cap G_{n_k}) + \mathbb{E}Y_{n_k} \ind(F_{n_k} \cap G^c_{n_k}) \nonumber\\
&\leq& (\mu_{n_k} + \epsilon)\mathbb{P}(F_{n_k}) + \mathbb{E}Y_{n_k} \ind(F_{n_k} \cap G^c_{n_k}). \label{up22}
\end{eqnarray}
The final term in~(\ref{up22}) is bounded above using the Cauchy-Schwarz inequality as
\begin{eqnarray}
\left(\mathbb{E}Y^2_{n_k}\right)^{\frac{1}{2}}\left(\mathbb{P}\left(F_{n_k} \cap G^c_{n_k}\right)\right)^{\frac{1}{2}}
&\leq& C\left(\mathbb{P}\left(F_{n_k} \cap G^c_{n_k}\right)\right)^{\frac{1}{2}} \nonumber\\
&\leq& C\left(\mathbb{P}\left(G^c_{n_k}\right)\right)^{\frac{1}{2}} \nonumber\\
&\leq& C\sqrt{\epsilon} \label{up_33}
\end{eqnarray}
for all~\(k\) large. Here~\(C  = \sup_n \mathbb{E}Y^2_n < \infty\) is a constant and the final estimate follows using~(\ref{cond_yn}).

Using~(\ref{up_33}) into~(\ref{up22}) we have
\[\mathbb{E}Y_{n_k}\ind(F_{n_k}) \leq (\mu_{n_k} + \epsilon)\mathbb{P}(F_{n_k}) + C \sqrt{\epsilon}\] for all~\(k\) large. Using the above in~(\ref{up1}),
we have
\begin{eqnarray}
\mathbb{E}Y_{n_k} &\leq& \mu_{n_k} + \epsilon\mathbb{P}(F_{n_k}) - \epsilon_0 \mathbb{P}(F^c_{n_k}) \nonumber\\
&\leq&\mu_{n_k} + \epsilon - \epsilon_0 \delta_0 \nonumber
\end{eqnarray}
for all~\(k\) large, where the final estimate follows using~(\ref{del_0_def}). This contradiction the definition that~\(\mu_n = \mathbb{E}Y_n.\)~\(\qed\)

\renewcommand{\theequation}{\thesection.\arabic{equation}}
\setcounter{equation}{0}
\section{Proof of Theorem~\ref{main_thm}}\label{pf1}
We first see that~\(T < \frac{1}{p_{av}}\) is not achievable. Suppose that~\(T < \frac{1}{p_{av}}\) is achievable.
We then have for any fixed~\(\epsilon  >0\) that
\begin{equation}\label{ach_t}
\mathbb{P}_{conf}\left(T(\underline{p}) \geq T+\epsilon\right) \leq \epsilon
\end{equation}
for all~\(n\) large. We therefore have
\begin{eqnarray}
\mathbb{E}_{conf}(T(\underline{p})) &=& \mathbb{E}_{conf}(T(\underline{p})) \ind(T(\underline{p}) < T+\epsilon)
 + \mathbb{E}_{conf}(T(\underline{p})) \ind(T(\underline{p}) \geq T+\epsilon) \nonumber\\
 &\leq& T + \epsilon  + (\mathbb{E}_{conf}T^2(\underline{p}))^{\frac{1}{2}} \mathbb{P}_{conf}\left(T(\underline{p}) \geq T+\epsilon\right)^{\frac{1}{2}} \nonumber\\
 &\leq& T + \epsilon  + C\sqrt{\epsilon} \label{ach_t2}
\end{eqnarray}
for some constant~\(C > 0.\) The final estimate is obtained from~(\ref{ach_t}) and the upper bound on the variance of~\(T(\underline{p})\) in~(\ref{var_tp_est_up}). Since~\(\epsilon  >0\) is arbitrary and~\(T < \frac{1}{p_{av}}\) this contradicts~(\ref{e_conf}).

We now show that that the~\(S(\underline{p})\) is arbitrarily close to~\(1\) if and only if~\(r(n) \longrightarrow 1\) as~\(n \rightarrow \infty.\) From~(\ref{prob_pi_p}), we have
\[\mathbb{P}^{(\pi,\underline{p})}(T_{det} < \infty)= 1-q_{\pi(1)}\ldots q_{\pi(r)}\] and so
\begin{eqnarray}
S(\underline{p}) &=& 1-\mathbb{E}_{sch}(q_{\pi(1)}\ldots q_{\pi(r)}) \nonumber\\
&\geq& 1-\frac{1}{r}\sum_{k=1}^{r}\mathbb{E}_{sch}q_{\pi(k)}^{r} \label{mid_main}
\end{eqnarray}
for all~\(n\) large. The middle inequality follows using the arithmetic-geometric inequality~\[x_1\ldots x_r \leq \frac{1}{r}\sum_{k=1}^{r}x_k^r\] for positive numbers~\(\{x_i\}.\)

Taking average over all configurations~\(\underline{p}\) we have
\begin{eqnarray}
\sum_{\underline{p}} S(\underline{p})\mathbb{P}_{conf}(\underline{p})
&\geq& 1-\mathbb{E}_{conf}\frac{1}{r}\sum_{k=1}^{r}\mathbb{E}_{sch}\left(q_{\pi(k)}^{r}\right) \label{mid_main2}\\
&=& 1-\frac{1}{r}\sum_{k=1}^{r}\mathbb{E}_{sch}\mathbb{E}_{conf}\left(q_{\pi(k)}^{r}\right) \nonumber\\
&=&1 - \mathbb{E}_{conf}(q_{1}^{r}) \nonumber\\
&\longrightarrow& 1 \label{fin_2main}
\end{eqnarray}
as~\(n \rightarrow \infty.\) The final estimate follows since~\(r(n) \longrightarrow \infty\) and~\(q_1  = 1-p_1 < 1\) since~\(p_1 > 0\) for all~\(p_1\) in the finite set~\(\Omega_{conf}.\)

Letting
\begin{equation}
A(\epsilon) := \{\underline{p}  : S(\underline{p}) > 1-\epsilon\},
\end{equation}
we evaluate
\begin{eqnarray}\label{two_sum}
\sum_{\underline{p}} S(\underline{p})\mathbb{P}_{conf}(\underline{p}) &=& I_1 + I_2
\end{eqnarray}
where
\[I_1 = \sum_{ \underline{p} \in A(\epsilon)} S(\underline{p})\mathbb{P}_{conf}(\underline{p}) \leq \mathbb{P}_{conf}(A(\epsilon))\] and
\[I_2 = \sum_{\underline{p} \notin A(\epsilon)} S(\underline{p}) \mathbb{P}_{conf}(\underline{p})  \leq \epsilon.\]

In particular, we have from~(\ref{two_sum}) and~(\ref{fin_2main}) that
\begin{equation}\label{a_eps_est}
\mathbb{P}_{conf}(A(\epsilon)) \geq I_1 \geq 1-2\epsilon
\end{equation}
for all~\(n\) large.

We now show that~\(\frac{1}{p_{av}}\) is achievable if and only if~\((i)\) and~\((ii)\) stated in Theorem~\ref{main_thm} hold.
Using Lemmas~\ref{lem1} and~\ref{lem2} and~(\ref{a_eps_est}) above, we have that if~\((i)-(ii)\) hold and~\(r(n) \longrightarrow \infty\) as~\(n \rightarrow \infty,\) then~\(\frac{1}{p_{av}}\) is achievable.

Suppose now that~\(\frac{1}{p_{av}}\) is achievable. For a fixed~\(\epsilon > 0,\) we have using~(\ref{cap_def}) that
\begin{equation}\label{ach_pp}
\mathbb{P}_{conf}\left(T(\underline{p}) \leq \frac{1}{p_{av}} + \epsilon\right) \longrightarrow 1
\end{equation}
as~\(n \rightarrow \infty.\) Using~(\ref{e_conf}), we obtain that
\begin{equation}\label{t_yn}
\mathbb{P}_{conf}\left(T(\underline{p}) \leq \mathbb{E}_{conf}T(\underline{p}) + 2\epsilon\right) \longrightarrow 1
\end{equation}
Using Lemma~\ref{lem_yn} with~\(Y_n = T(\underline{p}) = T_n(\underline{p})\) we have for every~\(\epsilon  >0\) that
\begin{equation}\label{t_yn}
\mathbb{P}_{conf}\left(T(\underline{p}) \geq \mathbb{E}_{conf}T(\underline{p}) - 2\epsilon\right) \longrightarrow 1
\end{equation}
as~\(n \rightarrow \infty.\) The Lemma~\ref{lem_yn} is applicable since~\(\sup_n \mathbb{E}_{conf}T^2(\underline{p}) < \infty\) using the upper bound in~(\ref{var_tp_est_up}) and the upper bound in~(\ref{e_conf}).

From the above we have that
\begin{equation}\label{conv_prob_t}
T(\underline{p}) - \mathbb{E}_{conf}T(\underline{p}) \longrightarrow 0
\end{equation}
in probability. From Lemma~\ref{lem2} we have that condition~\((ii)\) holds. Again using~(\ref{ach_pp}),~(\ref{conv_prob_t}) and~(\ref{e_conf}), we have that~\[\mathbb{E}_{conf}T(\underline{p}) \longrightarrow \frac{1}{p_{av}}\] as~\(n \rightarrow \infty.\) This implies from Lemma~\ref{lem1} that condition~\((i)\) holds.~\(\qed\)

\subsection*{Acknowledgement}
I thank Professors Rahul Roy and Federico Camia for crucial comments and for my fellowships.

\bibliographystyle{plain}

\end{document}